\documentclass{amsart}
\DeclareMathSymbol{\subsetneqq}     {\mathrel}{AMSb}{"24}
\DeclareMathSymbol{\subsetneq}      {\mathrel}{AMSb}{"28}

\def\Q{{\mathbf Q}}
\def\Z{{\mathbf Z}}

\def\R{{\mathbf R}}
\def\F{{\mathbf F}}

\def\Gal{\mathrm{Gal}}
\def\End{\mathrm{End}}
\def\Aut{\mathrm{Aut}}
\def\Hom{\mathrm{Hom}}
\def\I{{\mathcal I}}
\def\II{{{\mathcal I}_{v,X}}}
\def\J{{\mathcal J}}
\def\s{{\mathfrak S}}

\def\SL{\mathrm{SL}}

\def\M{\mathrm{M}}

\def\R{R}
\def\dim{\mathrm{dim}}
\def\O{{\mathcal O}}

\def\bmu{\boldsymbol \mu}
\newtheorem{thm}{Theorem}[section]
\newtheorem{lem}[thm]{Lemma}
\newtheorem{cor}[thm]{Corollary}

\theoremstyle{definition}
\newtheorem{defn}[thm]{Definition}
\newtheorem{ex}[thm]{Example}
\newtheorem{rem}[thm]{Remark}

\title[Semistable reduction of abelian varieties]
{Semistable reduction of abelian varieties over extensions
of small degree}

\author[A.\ Silverberg]{A.\ Silverberg*}
\thanks{* Partially supported by the National 
Security Agency.}
\address{Department of Mathematics, Ohio State University, 
231 W.\ 18 Avenue,
Columbus, Ohio 43210--1174, USA}
\email{silver\char`\@math.ohio-state.edu}

\author[Yu. G. Zarhin]{Yu. G. Zarhin**}
\thanks{** Partially supported by the National 
Science Foundation.}
\address{Department of Mathematics, Pennsylvania State University, 
University Park, PA 16802, USA,
\newline
\indent Institute for Mathematical Problems in Biology, 
Russian Academy of Sciences, Push\-chino, Moscow Region, 142292, Russia}
\email{zarhin\char`\@math.psu.edu}

\begin{document}

\begin{abstract}
We obtain necessary and sufficient conditions for abelian
varieties to acquire semistable reduction over fields of
low degree. Our criteria are expressed in terms of torsion
points of small order defined over unramified extensions.
\end{abstract}

\maketitle

\section{Introduction}

In this paper we obtain criteria for 
abelian varieties to acquire semistable reduction over fields
of certain given (small) degrees. Our criteria are expressed
in terms of unramified torsion points. 

Suppose that $X$ is an abelian variety defined over a field 
$F$, and $n$ is a positive integer not divisible by the characteristic
of $F$. Let $X^\ast$ denote the dual abelian variety of $X$, let
$X_n$ denote the kernel of multiplication by $n$ in $X(F^s)$, 
where $F^s$ denotes a separable closure of $F$, 
let
$X_n^\ast$ denote the kernel of multiplication by $n$ in $X^\ast(F^s)$,
and let $\bmu _n$ denote the $\Gal (F^s/F)$-module
of $n$-th roots of unity in $F^s$. 
The Weil pairing 
$e_n : X_n \times X_n^\ast \to {\bmu}_n$
is a $\Gal(F^s/F)$-equivariant nondegenerate pairing. 
If $S$ is a subgroup of $X_n$, let 
$$S^{\perp_n} = 
\{ y \in X_n^\ast : e_n(x,y) = 1 \text{ for every } x \in S \}
\subseteq X_n^\ast.$$
For example, if $n = m^2$ and $S = X_m$, then $S^{\perp_n} = X_m^\ast$.
If $X$ is an elliptic curve and $S$ is a cyclic subgroup of
order $n$, then $S^{\perp_n} = S$.
Suppose that $v$ is a discrete valuation on $F$ whose residue
characteristic does not divide $n$. 

Previously we showed that if $n \ge 5$
then $X$ has semistable reduction at $v$ if and only if
there exists a subgroup $S$ of $X_n$ such that all the points on
$S$ and on $S^{\perp_n}$ are defined over an extension of $F$
unramified over $v$ (see Theorem 4.5 of \cite{dpp}; 
see also Theorem 6.2 of \cite{semistab}). 
In the current paper we show that if 
there exists a subgroup $S$ of $X_n$, 
for $n = 2$, $3$, or $4$ (respectively),
such that all the points on
$S$ and on $S^{\perp_n}$ are defined over an extension of $F$
unramified over $v$, then $X$ acquires semistable reduction
over every degree $4$, $3$, or $2$ (respectively)
extension of $F$ totally ramified above $v$.
We also give necessary and sufficient conditions for 
semistable reduction over quartic, cubic, and quadratic extensions. 
Namely, if $L$ is a totally ramified extension of $F$ 
of degree $4$, $3$, or $2$, respectively,
then $X$ has semistable reduction over $L$
if and only if there exist a finite unramified
extension $K$ of $F$, an abelian variety $Y$ over $K$ which
is $K$-isogenous to $X$, and a subgroup $S$ of $Y_n$, 
for $n = 2$, $3$, or $4$, respectively, 
such that all the points of $S$ and of $S^{\perp_n}$ are
defined over an unramified extension of $K$.
If $X$ is an elliptic curve one may take $Y = X$. This
is not true already for abelian surfaces. However, one may take 
$Y = X$ in the special case where $X$ has purely additive and
potentially good reduction, with no restriction on the dimension. 

The study of torsion subgroups of 
abelian varieties with purely additive reduction was 
initiated in \cite{LenstraOort}
and pursued in \cite{Lorenzini} (see \cite{Frey} and  
\cite{Flexor-Oesterle} for the case of elliptic curves).
See \cite{Kraus} for a study of the smallest extension over
which an elliptic curve with additive and potentially good 
reduction acquires good reduction.

We state and prove Theorem \ref{oneway} in
the generality $n \ge 2$ (rather than just $2 \le n \le 4$) 
since doing so requires no extra work
and affords us the opportunity to give 
a slightly different exposition from 
that in \cite{dpp} for $n \ge 5$, which highlights the
method. See \S\ref{ssredsect} for
our major results, see \S\ref{appsect} for applications and
refinements, and see \S\ref{exssect} for examples which
demonstrate that our results are sharp.

\section{Notation and definitions}

Define 
$$R(n) = 1 \text{ if } n \ge 5, \quad R(4) = 2, 
\quad R(3) = 3, \quad R(2) = 4.$$

If $X$ is an abelian variety over a field $F$, and 
$\ell$ is a prime not equal to the characteristic of $F$,
let
$$\rho_{\ell,X} : \Gal(F^s/F) \to \Aut(T_\ell(X))$$
denote the $\ell$-adic representation on the Tate module $T_\ell(X)$
of $X$. We will write $\rho_\ell$ when there is no ambiguity.
Let $V_\ell(X) = T_\ell(X)\otimes_{\Z_\ell}\Q_\ell$.

If $L$ is a Galois extension of $F$
and $w$ is an extension of $v$ to $L$,
let $\I(w/v)$ denote the inertia subgroup at $w$ of $\Gal(L/F)$.
Throughout this paper we will let $\I$ denote 
$\I({\bar v}/v)$, where ${\bar v}$ is a fixed extension of $v$ to $F^s$,
and we will let $\J$ denote the first ramification group 
(i.e., the wild inertia group). We also write $\I_{w}$ for
$\I({\bar v}/w)$.

\begin{defn}
Suppose $L/F$ is an extension of fields, $w$ is a discrete valuation
on $L$, and $v$ is the restriction of $w$ to $F$. 
Let $e(w/v) = [w(L^\times):v(F^\times)]$. We say that
$w/v$ is {\em unramified} if $e(w/v) = 1$ 
and the residue field extension is separable.
We say that
$w/v$ is {\em  totally ramified} if $w$ is the unique extension
of $v$ to $L$ and the residue field extension is purely inseparable.
We say that
$w/v$ is {\em tamely ramified} if the residue field extension 
is separable and $e(w/v)$ is not divisible by the residue characteristic.
\end{defn}

\section{Preliminaries}

\begin{thm}
\label{quasithm}
Suppose $n$ is an integer, $n \ge 2$, $\O$ is an integral 
domain of characteristic zero such that no rational prime which 
divides $n$ 
is a unit in $\O$, $\alpha \in \O$, $\alpha$ has finite 
multiplicative order, 
and $(\alpha-1)^2 \in n\O$.
Then $\alpha^{R(n)} = 1$.
\end{thm}

\begin{proof}
See Corollary 3.3 of \cite{serrelem}.
\end{proof}

\begin{lem}[Lemma 5.2 of \cite{semistab}]
\label{localglobal}
Suppose that $d$ and $n$ are positive integers, and for each prime $\ell$ 
which divides
$n$ we have a matrix $A_\ell \in M_{2d}(\Z_\ell)$ such that the
characteristic polynomials of the $A_\ell$ have integral coefficients
independent of $\ell$, and such that $(A_\ell-1)^2 \in nM_{2d}(\Z_\ell)$.
Then for every eigenvalue $\alpha$ of $A_\ell$, $(\alpha-1)/\sqrt{n}$
satisfies a monic polynomial with integer coefficients.
\end{lem}

\begin{thm}[Galois Criterion for Semistable Reduction]
\label{galcrit}
Suppose $X$ is an abel\-ian variety over a field $F$, $v$ is a discrete 
valuation on $F$, and $\ell$ is a prime not equal to the residue
characteristic of $v$. 
Then the following are equivalent:
\begin{enumerate}
\item[(i)] $X$ has semistable reduction at $v$,
\item[(ii)] $\I$ acts unipotently on $T_\ell(X)$; i.e.,
all the eigenvalues of $\rho_{\ell}(\sigma)$ are $1$, 
for every $\sigma \in \I$,
\item[(iii)] for every $\sigma \in \I$, 
$(\rho_{\ell}(\sigma)-1)^2 = 0$.
\end{enumerate}
\end{thm}

\begin{proof}
See Proposition 3.5 and Corollaire 3.8 of \cite{SGA} 
and Theorem 6 on p.~184 of \cite{BLR}.
\end{proof}

\begin{lem}
\label{primrtsof1}
Suppose $\ell$ is a prime number and $\zeta$ is a primitive 
$\ell^s$-th root of unity. Then 
$$\frac{(\zeta-1)^{\varphi(\ell^s)}}{\ell}$$
is a unit in $\Z[\zeta]$.
\end{lem}

\begin{proof}
See for example the last two lines on p.~9 of \cite{Wash}.
\end{proof}

\section{Lemmas}
\label{ssredlemsect}

\begin{rem}
\label{cyclicrem}
Suppose $w$ is a discrete valuation on a field $L$, $L$ is a finite
extension of a field $F$, $v$ is the restriction of $w$ to $F$, 
and $w/v$ is totally and tamely ramified. 
Then the maximal unramified extension $L_{nr}$ of $L$ is the compositum
of $L$ with the maximal unramified extension $F_{nr}$ of $F$.
Further, $L_{nr}/F_{nr}$ is a cyclic extension whose degree is
$[L:F]$ (see \S8 of \cite{Frohlich}, especially Corollary 3 on p.~31).
Since passing to the maximal unramified extensions does not change
the inertia groups, 
it follows that $\I_{w}$ is a normal subgroup of $\I$, and 
$\I/\I_{w}$ is cyclic of order $[L:F]$.
\end{rem}

\begin{lem}
\label{rootsofone}
Suppose $v$ is a discrete valuation on a field $F$ with residue
characteristic $p \ge 0$, $R$ is a positive integer, $\ell$ is a prime,
$p$ does not divide $R\ell$,
and $L$ is a degree $R$ extension of $F$
which is totally ramified above $v$. 
Suppose that $X$ is an abelian variety over $F$,
and for every $\sigma \in \I$,
all the eigenvalues of $\rho_\ell(\sigma)$ are $R$-th roots of unity.
 Then $X$ has
semistable reduction at the extension of $v$ to $L$.
\end{lem}

\begin{proof}
This was proved in Lemma 5.5 of \cite{semistab} in the case where
$L$ is Galois over $F$. However, the same proof also works in general.
This follows from the fact that in the proof we replaced $F$ by its
maximal unramified extension. For fields which have no non-trivial
unramified extensions, every totally and tamely ramified extension is 
cyclic (and therefore Galois), and for each degree prime to the residue 
characteristic,
there is a unique totally ramified extension of that degree.
See \S8 of \cite{Frohlich}, especially Corollary 3 on p.~31.
\end{proof}

The following result yields a converse of Theorem 5.1 of \cite{serrelem}.

\begin{lem}
\label{algprop}
Suppose $\O$ is an integral domain of characteristic zero, and
$\ell$ is a prime number. Suppose $k$, $r$, and $m$ are positive integers
such that $k \ge m\varphi(\ell^r)$. Suppose $\alpha \in \O$ and
$\alpha^{\ell^r} = 1$. Then 
$(\alpha-1)^{k} \in \ell^m\Z[\alpha]$.
\end{lem}

\begin{proof}
Let $s$ be the smallest positive integer such that $\alpha^{\ell^s} = 1$. 
Then  
$$(\alpha-1)^{k} \in (\alpha-1)^{m\varphi(\ell^s)}\Z[\alpha] 
\subseteq \ell^m\Z[\alpha],$$
by Lemma \ref{primrtsof1}.
\end{proof}

\begin{lem}
\label{ssprelem}
Suppose $X$ is an abelian variety over a field $F$, $v$ is a discrete 
valuation on $F$, $n$ and $m$ are integers, and $n$ is not 
divisible by the residue characteristic of $v$. 
Suppose $\sigma \in \I$.
If there exists a subgroup $S$ of $X_n$
such that $(\sigma^m-1)S = 0$ and $(\sigma^m-1)S^{\perp_n} = 0$,
then 
$(\sigma^m-1)^2 X_n = 0$.
\end{lem}

\begin{proof}
The map $x \mapsto (y \mapsto e_n(x,y))$ induces a
$\Gal(F^s/F)$-equivariant isomorphism from $X_n/S$ onto
$\Hom(S^{\perp_n},\bmu_n)$. Since  
$\sigma = 1$ on $\bmu_n$, and 
$\sigma^m = 1$ on $S^{\perp_n}$, it follows that $\sigma^m = 1$ on $X_n/S$. 
Therefore, 
$(\sigma^m-1)^2X_n \subseteq (\sigma^m-1)S = 0$.
\end{proof}

\begin{lem}
\label{sslem}
Suppose $X$ is an abelian variety over a field $F$, $v$ is a discrete 
valuation on $F$, $n$ is an integer not 
divisible by the residue characteristic of $v$,  
and $S = X_n^\I$.  
Then $\I$ acts as the identity on $S^{\perp_n}$
if and only if 
$(\sigma-1)^2 X_n = 0$ for every $\sigma \in \I$.
\end{lem}

\begin{proof}
Applying Lemma \ref{ssprelem} with $m = 1$, we obtain the
forward implication. 

Conversely, suppose that
$(\sigma-1)^2 X_n = 0$ for every $\sigma \in \I$. Writing
$\sigma^n = ((\sigma-1)+1)^n$, it is easy to see that 
$\sigma^n = 1$ on $X_n$ for every $\sigma \in \I$. 
Since $n$ is not divisible by the residue characteristic of $v$, 
$X_n$ and $X_n^\ast$ are tamely ramified at $v$.
Then the action of $\I$ on $X_n$ and on $X_n^\ast$ factors through 
the tame inertia group $\I/\J$. 
Let $\tau$ denote
a lift to $\I$ of a topological generator of the pro-cyclic group
$\I/\J$. Since 
$$e_n((\tau-1)X_n,(X_n^\ast)^\I) = 1,$$
we have 
$$\#((X_n^\ast)^\I)\#((\tau-1)X_n) \le \#X_n^\ast.$$
The map from $X_n$ to $(\tau-1)X_n$ defined by
$y \mapsto (\tau-1)y$ defines a short exact sequence
$$0 \to S \to X_n \to (\tau-1)X_n \to 0.$$
Therefore, 
$$\#S\#((\tau-1)X_n) = \#X_n = \#S\#S^{\perp_n}.$$
Similarly, 
$$\#((X_n^\ast)^\I)\#((\tau-1)X_n^\ast) = \#X_n^\ast.$$
Therefore, 
$$\#S^{\perp_n} = \#((\tau-1)X_n) \le \#((\tau-1)X_n^\ast).$$
Since $(\tau-1)X_n^\ast \subseteq S^{\perp_n}$,
we conclude that 
$$S^{\perp_n} = (\tau-1)X_n^\ast.$$ From the natural
$\Gal(F^s/F)$-equivariant isomorphism 
$X_n^\ast \cong \Hom(X_n,\bmu_n)$ it follows
that $(\tau-1)^2X_n^\ast = 0$.
Therefore, $\I$ acts as the identity on $S^{\perp_n}$.
\end{proof}

\begin{lem}
\label{sslem2}
Suppose $X$ is an abelian variety over a field $F$, $v$ is a discrete 
valuation on $F$, and $n$ is an integer not 
divisible by the residue characteristic of $v$. 
If $X$ has semistable reduction at $v$, then
\begin{enumerate}
\item[(i)] $(\sigma-1)^2X_n = 0$ for every $\sigma \in \I$,
\item[(ii)] $\I$ acts as the identity on $(X_n^\I)^{\perp_n}$,
\item[(iii)] $(\sigma^n-1)X_n = 0$ for every $\sigma \in \I$;
in particular, $X_n$ is tamely ramified at $v$.
\end{enumerate}
\end{lem}

\begin{proof}
By Theorem \ref{galcrit}, we have (i). By Lemma \ref{sslem}, we have (ii).
In the proof of Lemma \ref{sslem}, we showed that (i) implies (iii).
\end{proof}

\begin{lem}
\label{tame}
Suppose $X$ is an abelian variety over a field $F$, $v$ is a discrete 
valuation on $F$ of residue characteristic $p \ge 0$, and 
$\ell$ is a prime number not equal to $p$. 
If $X_\ell$ is tamely ramified at $v$, then $T_\ell(X)$ is
tamely ramified at $v$.
\end{lem}

\begin{proof}
If $p = 0$ then the wild inertia group $\J$ is 
trivial and we are done.
Suppose $p > 0$ and $\sigma \in \J$. 
Since $p \ne \ell$, $\rho_{\ell}(\J)$ is a finite $p$-group.
Therefore, $\rho_{\ell}(\sigma)$ has order a power of $p$.
Since $X_\ell$ is tamely ramified, 
$\rho_{\ell}(\sigma)-1 \in \ell\End(T_\ell(X))$.
It follows that $\rho_{\ell}(\sigma) = 1$ if $\ell \ge 3$,
and $\rho_{\ell}(\sigma)^2 = 1$ if $\ell = 2$. 
Since $p$ and $\ell$ are relatively prime, 
$\rho_{\ell}(\sigma) = 1$.
\end{proof}

\begin{lem}
\label{lemT}
Suppose $X$ is an abelian variety over a field $F$, $n = 2$, $3$, or 
$4$, $\ell$ is the prime divisor of $n$,
$v$ is a discrete valuation on $F$ whose residue characteristic is not
$\ell$, $t$ is a non-negative integer,
$L$ is an extension of $F$ of degree $R(n)^{t+1}$
which is totally ramified above $v$, 
and $X$ has semistable reduction
over $L$ above $v$. Let $\tau$ denote a lift to $\I$ of a topological 
generator of the pro-cyclic group $\I/\J$. 
Let $\gamma = \rho_{\ell}(\tau)^{\R(n)^t}$, 
let $\lambda = (\gamma - 1)^2/n$, and let
$$T = T_\ell(X) + {\lambda}T_\ell(X) + {\lambda^2}T_\ell(X) + \cdots + 
{\lambda^{R(n)-1}}T_\ell(X).$$
Then:
\begin{enumerate}
\item[(a)]  $T$ is the smallest $\lambda$-stable $\Z_\ell$-lattice 
in $V_\ell(X)$ which contains $T_{\ell}(X)$,
 \item[(b)]  $(\gamma^{R(n)} - 1)^{2} = 0$,
 \item[(c)]  $n^{R(n)-1}T \subseteq T_{\ell}(X) \subseteq T$,
 \item[(d)]  $(\gamma - 1)^{2R(n)} \subseteq nT_{\ell}(X)$,
 \item[(e)]  if $n = 2$ or $3$, then 
 $nT \subseteq T_{\ell}(X)$ if and only if 
 $(\gamma - 1)^{4}T_{\ell}(X) \subseteq nT_{\ell}(X)$,
 \item[(f)]  if $n = 2$, then $4T \subseteq T_{2}(X)$ if and only if
 $(\gamma - 1)^{6}T_{2}(X) \subseteq 2T_{2}(X)$,
 \item[(g)]  if $n = 4$, then $2T \subseteq T_{2}(X)$ if and only if
 $(\gamma - 1)^{2}T_{2}(X) \subseteq 2T_{2}(X)$. 
\end{enumerate}
\end{lem}

\begin{proof}
Let $w$ denote the restriction of ${\bar v}$ to $L$.
By Remark \ref{cyclicrem}, $\I/\I_{w}$ is cyclic of
order $R(n)^{t+1}$. 
By Theorem \ref{galcrit}, we have (b).
It follows that
$(\lambda + \gamma)^2(\lambda + \gamma - 1)^2 = 0$ if $n = 2$, 
$\lambda(\lambda + \gamma)^2 = 0$ if $n = 3$, 
and
$\lambda(\lambda + \gamma) = 0$ if $n = 4$. 
Therefore, $\lambda$ satisfies a polynomial over
$\Z[\gamma]$ of degree $R(n)$, and we have (a) and (c). 
From the definition of $T$ we easily deduce (e), (f), and (g).
Further, (d) follows from (b).
\end{proof}

We will apply the following result only in Corollary \ref{4326cor}e.

\begin{thm}
\label{divby23prop}
Suppose $L/F$ is a finite separable field extension, 
$w$ is a discrete valuation on $L$, 
and $v$ is the restriction of $w$ to $F$. 
Suppose $X$ is a $d$-dimensional abelian variety over $F$
which has semistable reduction at $w$ but not at $v$.
Then $[\I_{v}:\I_{w}]$ has a prime divisor $q$ such that 
$q \le 2d + 1$.
\end{thm}

\begin{proof}
Let $\ell$ be a prime not equal to the residue characteristic $p$,
and let 
$$\II = 
\{\sigma \in \I_{v} : 
\sigma \text{ acts unipotently on } V_{\ell}(X) \}.$$
We have  
$\I_{w} \subseteq \II \subsetneqq \I_{v}$ by Theorem \ref{galcrit}, 
since $X$ has semistable reduction at $w$ but not at $v$. 
Let $F_{v}$ be the completion of $F$ at $v$ and let
$F_{v}^{nr}$ be the maximal unramified extension 
of $F_{v}$.  
Then $\II$ is an open normal subgroup of $\I_{v}$, 
is independent of $\ell$,  
and cuts out the smallest Galois extension $F'$ of $F_{v}^{nr}$ 
over which $X$ has semistable reduction (see pp.~354--355
of \cite{SGA}). We have $\Gal(F'/F_{v}^{nr}) \cong \I_{v}/\II$. 
By a theorem of Raynaud (see Proposition 4.7 of \cite{SGA}),
$X$ has semistable reduction over $F_{v}^{nr}(X_{n})$,
for every integer $n$ not divisible by $p$ and greater than $2$.
The intersection $M$ of these fields therefore contains $F'$.
As on the top of p.~498 of \cite{SerreTate}, 
every prime divisor of $[M:F_{v}^{nr}]$ is at most $2d+1$
(see Theorem 4.1 and Formula 3.1 of \cite{JPAA} for 
an explicit integer that $[M:F_{v}^{nr}]$ divides).
Thus, if $q$ is a prime divisor of $[\I_{v}:\II]$ 
then $q \le 2d+1$. 
Since $\I_{w} \subseteq \II \subsetneqq \I_{v}$, 
we obtain the desired result.
\end{proof}

\begin{rem}
\label{divby23}
With hypotheses and notation as in Theorem \ref{divby23prop},
let $k_{w}$ and $k_{v}$ denote the residue fields.
Then $[\I_{v}:\I_{w}] = e(w/v)[k_{w}:k_{v}]_{i}$, where
the subscript $i$ denotes the inseparable degree (see
Proposition 21 on p.~32 of \cite{Corps} for the case where
$L/F$ is Galois. In the non-Galois case, take a Galois extension
$L'$ of $F$ which contains $L$, and apply the result to
$L'/L$ and $L'/F$, to obtain the result for $L/F$).
Taking completions, then $[L_{w}:F_{v}] = e(w/v)[k_{w}:k_{v}]
= [\I_{v}:\I_{w}][k_{w}:k_{v}]_{s}$, where
the subscript $s$ denotes the separable degree.
Therefore, the prime $q$ from Theorem \ref{divby23prop}
divides $[L_{w}:F_{v}]$.
\end{rem}

\section{Semistable reduction}
\label{ssredsect}

The results in this section extend the results of \cite{dpp}
to the cases $n = 2, 3, 4$. Theorem \ref{oneway} is also
a generalization of Corollary 7.1 of \cite{semistab}.

\begin{rem}
\label{lemrem}
Suppose $X$ is an abelian variety over a field $F$, $v$ is a discrete 
valuation on $F$, and $n$ is an integer greater than $1$ which is not 
divisible by the residue characteristic of $v$. 
By Lemma \ref{sslem}, the following two statements are equivalent:
\begin{enumerate}
\item[(a)] 
there exists a subgroup $S$ of $X_n$
such that $\I$ acts as the identity on $S$ and on $S^{\perp_n}$,
\item[(b)] $(\sigma - 1)^2X_n = 0$ for every $\sigma \in \I$.
\end{enumerate}
\end{rem}

\begin{thm}
\label{oneway}
Suppose $X$ is an abelian variety over a field $F$, $v$ is a discrete 
valuation on $F$, and $n$ is an integer greater than $1$ which is not 
divisible by the residue characteristic of $v$. 
Suppose there exists a subgroup $S$ of $X_n$
such that $\I$ acts as the identity on $S$ and on $S^{\perp_n}$.
Then $X$ has semistable reduction over every degree $R(n)$
extension of $F$ totally ramified above $v$.
\end{thm}

\begin{proof}
Suppose $\sigma \in \I$.
By Lemma \ref{sslem}, $(\sigma-1)^2 X_n = 0$.
Let $\I' \subseteq \I$ be the inertia group for the prime 
below ${\bar v}$ in a
finite Galois extension of $F$ over which $X$ has 
semistable reduction.
Then $\sigma^r \in \I'$ for some $r$.
Let $\ell$ be a prime divisor of $n$. 
Theorem \ref{galcrit} implies that
$(\rho_{\ell}(\sigma)^r-1)^2 = 0$.
Let $\alpha$ be an eigenvalue of $\rho_{\ell}(\sigma)$.
Then $(\alpha^r-1)^2 = 0$.
Therefore, $\alpha^r = 1$.
By our hypothesis, 
$$(\rho_{\ell}(\sigma)-1)^2 \in n\M_{2d}(\Z_\ell),$$
where $d = \dim(X)$.
By Th\'eor\`eme 4.3 of \cite{SGA}, 
the characteristic polynomial of $\rho_{\ell}(\sigma)$
has integer coefficients which are independent of $\ell$.
By Lemma \ref{localglobal},
$(\alpha-1)^2 \in n{\bar \Z}$, where ${\bar \Z}$ denotes
the ring of algebraic integers.
By Theorem \ref{quasithm} we have $\alpha^{R(n)} = 1$.
The result now follows from Lemma \ref{rootsofone}.
\end{proof}

\begin{cor}[Theorem 4.5 of \cite{dpp}]
\label{fromdpp}
Suppose $X$ is an abelian variety over a field $F$, $v$ is a discrete 
valuation on $F$, $n$ is an integer not 
divisible by the residue characteristic of $v$, and $n \ge 5$. 
Then $X$ has semistable reduction at $v$ if and only if
there exists a subgroup $S$ of $X_n$
such that $\I$ acts as the identity on $S$ and on $S^{\perp_n}$.
\end{cor}

\begin{proof}
If $X$ has semistable reduction at $v$, then by Theorem \ref{galcrit},
$(\sigma - 1)^2X_n = 0$ for every $\sigma \in \I$. Apply 
Lemma \ref{sslem}.

For the converse, apply Theorem \ref{oneway} with $n \ge 5$.
\end{proof}

\begin{rem}
\label{otherway}
It follows immediately from Theorem \ref{galcrit} and Lemma \ref{sslem}
that if $X$ has semistable reduction above $v$ over a degree 
$m$ extension of $F$ totally ramified above $v$,
then there exists a subgroup $S$ of $X_n$
such that $\I$ acts via a cyclic quotient of order $m$ 
on $S$ and on $S^{\perp_n}$. (If $L$ is the extension of $F$,
let $w$ be the restriction of ${\bar v}$
to $L$ and let $S = X_n^{\I_w}$.) 
Theorem \ref{bothways} below gives a different result in the
direction converse to Theorem \ref{oneway}, and, further, gives
conditions for semistable reduction which are both necessary
and sufficient, thereby giving a generalization of 
Corollary \ref{fromdpp} to the cases $n = 2, 3, 4$.
Note that in the case $n \ge 5$, the equivalence of (i)
and (ii) in Theorem \ref{bothways} is just a restatement
of Corollary \ref{fromdpp} (since $R(n) = 1$ if $n \ge 5$).
We remark that in that case, one can take
(in the notation of Theorem \ref{bothways}) $Y = X$ and 
$\varphi$ the identity map.
\end{rem}

\begin{thm}
\label{bothways}
Suppose $n = 2$, $3$, or $4$, respectively. Suppose 
$X$ is an abelian variety over a field $F$, and 
$v$ is a discrete valuation on $F$ whose residue characteristic 
does not divide $n$. 
Suppose $t$ is a non-negative integer and
$L$ is an extension of $F$ of degree $\R(n)^{t+1}$ which 
is totally ramified above $v$. 
Then the following are equivalent:
\begin{enumerate}
\item[(i)] $X$ has semistable reduction over $L$ above $v$,
\item[(ii)] there exist an abelian variety $Y$ over a finite
extension $K$ of $F$ unramified above $v$, a separable $K$-isogeny
$\varphi : X \to Y$, and a subgroup 
$S$ of $Y_n$ such that 
$\I$ acts via a cyclic quotient of order $\R(n)^t$ on $S$ and 
on $S^{\perp_n}$.
\end{enumerate}
One can take $\varphi$ so that
its kernel is killed by $8$, $9$, or $4$, respectively.
If $X$ has potentially good reduction at $v$, then 
one can take $\varphi$ so that its kernel is killed by $2$,
$3$, or $2$, respectively.
\end{thm}

\begin{proof}
Let $\ell$ denote the prime divisor of $n$.

Suppose $K$ is a finite extension of $F$ unramified above $v$,
$Y$ is an abelian variety over $K$, $X$ and $Y$ are $K$-isogenous, and 
$S$ is a subgroup of $Y_n$ such that
$\I$ acts via a cyclic quotient of order $\R(n)^t$ on $S$ and 
on $S^{\perp_n}$. Suppose $\sigma \in \I$. 
By Lemma \ref{ssprelem}, $(\sigma^{\R(n)^t}-1)^2Y_n = 0$, i.e.,
$$(\rho_{\ell,Y}(\sigma^{\R(n)^t})-1)^2 \in nM_{2d}(\Z_\ell).$$
Let $\alpha$ be an eigenvalue of $\rho_{\ell,Y}(\sigma)$.
Since $Y$ has potentially semistable reduction, $\alpha$ is a root
of unity. 
By Theorem \ref{quasithm}, $(\alpha^{\R(n)^t})^{\R(n)} = 1$.
Therefore, all eigenvalues of $\rho_{\ell,Y}(\sigma)$ are 
${\R(n)}^{t+1}$-th
roots of unity. By Lemma \ref{rootsofone}, $Y$ has semistable reduction
over $LK$ above $v$. Since $X$ and $Y$ are $K$-isogenous and
$K/F$ is unramified above $v$, 
$X$ has semistable reduction over $L$ above $v$. 

Conversely, suppose $X$ has semistable reduction over $L$ above $v$.
By Lemma \ref{sslem2}iii, 
for every $\sigma \in \I$ we have
$(\sigma^{nR(n)^{t+1}}-1)X_n = 0$.
Since $nR(n)^{t+1}$ is 
not divisible by the residue characteristic, 
$X_n$ is tamely ramified at $v$. Then the
action of $\I$ on $X_n$ factors through $\I/\J$. 
Let $\tau$ denote a lift to $\I$ of a topological generator of the 
pro-cyclic group $\I/\J$. 
Let $T$ denote the $\Z_\ell$-lattice obtained from Lemma \ref{lemT}.
By Lemma \ref{tame}, $T$ is stable under $\I$. 
Note that $n^{R(n)-1} = 8$, $9$, or $4$ when
$n = 2$, $3$, or $4$, respectively.
Let $C = T/T_\ell(X)$, and 
view $C$ as a subgroup of 
$X_8$, $X_9$, or $X_4$, respectively. 
Let $Y = X/C$.
Then the projection map $X \to Y$ is a separable isogeny defined 
over a finite separable extension $K$ of $F$
which is unramified over $v$, 
$$T_\ell(Y) = T, \qquad \text{ and } \qquad 
(\rho_{\ell,Y}(\tau)^{\R(n)^t}-1)^2Y_n = 0.$$
Let $K'$ (respectively, $L'$) be the maximal unramified extension 
of $K$ (respectively, $L$) in $F^{s}$,  
let $M$ be the degree ${\R(n)^t}$ extension of $K'$ in $K'L'$ cut out by 
$\tau^{\R(n)^t}$, let $w$ be the restriction of ${\bar v}$ to $M$, 
and let $S = Y_n^{\I_w}$. 
Then $\tau^{\R(n)^t}$ is a lift to 
$\I_w$ of a topological generator of the pro-cyclic group
$\I_w/\J_w$, where $\J_w$ is the first ramification group of $\I_w$.
By Lemma \ref{sslem}, $\tau^{\R(n)^t}$ acts as the identity on
$S$ and on $S^{\perp_n}$. Therefore, $\I$ acts on $S$ and on 
$S^{\perp_n}$ via the cyclic group $\I/\I_w \cong \Gal(M/K')$.

As in Lemma \ref{lemT}, let $\gamma = \rho_{\ell,X}(\tau)^{\R(n)^t}$ 
and let $\lambda = (\gamma - 1)^2/n$.
If $X$ has potentially good reduction at $v$, then
$\gamma^{R(n)} = 1$. 
Let $\mu = \lambda + \gamma$.
Then $\mu^2 = \mu$ and $T = T_\ell(X) + \mu T_\ell(X)$.
Since $\mu = (\gamma^2 + 1)/2$ if $n = 2$, 
$\mu = (\gamma^2+\gamma+1)/3$ if $n = 3$, and 
$\mu = (\gamma + 1)/2$ if $n = 4$, 
it follows that 
$C$ is a subgroup of $X_2$, $X_3$, or $X_2$, respectively.
\end{proof}

Since the most interesting case of 
Theorem \ref{bothways} is the case $t = 0$, we explicitly 
state that case.

\begin{cor}
\label{bothcor}
Suppose $n = 2$, $3$, or $4$, respectively. Suppose  
$X$ is an abelian variety over a field $F$, and 
$v$ is a discrete valuation on $F$ whose residue characteristic 
does not divide $n$. 
Suppose
$L$ is an extension of $F$ of degree $4$, $3$, or
$2$, respectively, which 
is totally ramified above $v$. 
Then the following are equivalent:
\begin{enumerate}
\item[(i)] $X$ has semistable reduction over $L$ above $v$,
\item[(ii)] there exist an abelian variety $Y$ over a finite
extension $K$ of $F$ unramified above $v$, a separable
$K$-isogeny $\varphi : X \to Y$, 
and a subgroup 
$S$ of $Y_n$ such that 
$\I$ acts as the identity on $S$ and 
on $S^{\perp_n}$.
\end{enumerate}
Further, $\varphi$ can be taken so that
its kernel is killed by $8$, $9$, or $4$, respectively.
If $X$ has potentially good reduction at $v$, then $\varphi$ 
can be taken so that its kernel is killed by $2$, $3$, or
$2$, respectively.
\end{cor}

\section{Applications and refinements}
\label{appsect}

In the next result we show that the numbers in Theorem \ref{bothways}
and Corollary \ref{bothcor} can be improved for abelian varieties
of dimension $1$, $2$ (if $n = 2$ or $3$), and $3$ (if $n = 2$).
In \S\ref{exssect} we show that the numbers in Theorem \ref{ellcor} 
are sharp. See also \cite{Katz}, which deals with other 
problems concerned with finding a ``good'' abelian variety in an isogeny
class, with an answer depending on the dimension.

\begin{thm}
\label{ellcor}
In Theorem \ref{bothways} and Corollary \ref{bothcor},
with $d = \dim(X)$, 
$\varphi$ can be taken so that its kernel is killed
by $4$ if $d = 3$ and $n = 2$,
by $3$ if $d =2$ and $n = 3$, and 
by $2$ if $d = n = 2$. If $d = 1$, then we can take
$Y = X$ and $\varphi$ the identity map.
\end{thm}

\begin{proof}
We use the notation from Lemma \ref{lemT} and from the proof of 
Theorem \ref{bothways}.

Suppose $n = 2$ or $3$.
By Lemma \ref{lemT}d, $\gamma$ acts unipotently on the
$\F_{\ell}$-vector space 
$X_{\ell} \cong \frac{1}{\ell}T_{\ell}(X)/T_{\ell}(X)$.
Therefore,
$(\gamma - 1)^{2d}X_{\ell} = 0$. 
By Lemma \ref{lemT}e, if $d = 2$ then
$C$ is killed by $n$. By Lemma \ref{lemT}f, if $n = 2$ and $d = 3$, 
then $C$ is killed by $4$. 
If $d = 1$, then $\lambda$ is an endomorphism of $T_{\ell}(X)$,
so $T = T_{\ell}(X)$ and $Y = X$. 

Suppose $d = 1$ and $n = 4$. Since $\tau \in \I$, we have 
$\gamma \in \SL_{2}(\Z_{2})$. Therefore, the eigenvalues of 
$\gamma$ are either both $1$ or both $-1$. Therefore either
$(\gamma - 1)^{2} = 0$ or $(\gamma + 1)^{2} = 0$. In both cases,
$(\gamma - 1)^{2}X_{4} = 0$. Therefore, 
$\lambda$ is an endomorphism of $T_{2}(X)$ and 
$Y = X$.
\end{proof}

We can therefore take $Y = X$ in Theorem \ref{bothways} and 
Corollary \ref{bothcor} when $X$ is an elliptic curve.
This is not the case in general for abelian varieties
of higher dimension, as shown by the examples in the next section.
However, in Corollary \ref{paddcor} below we will show
that a result of this sort does hold for abelian varieties
with purely additive potentially good reduction. 

Next, we will give criteria for an elliptic curve to
acquire semistable reduction over extensions of degree
$2$, $3$, $4$, and either $6$ or $12$.

\begin{cor}
\label{4326cor}
Suppose 
$X$ is an elliptic curve over a field $F$, and 
$v$ is a discrete valuation on $F$ of residue characteristic 
$p \ge 0$. 
\begin{enumerate}
\item[(a)] If $p \ne 2$, 
then $X$ has semistable reduction above $v$ over a 
totally ramified quartic extension of $F$  
if and only if 
$X$ has an $\I$-invariant point of order $2$.
\item[(b)] If $p \ne 3$, 
then $X$ has semistable reduction above $v$ over a totally ramified  
cubic extension of $F$  
if and only if $X$ has an $\I$-invariant point of order $3$.
\item[(c)] 
If $p \ne 2$,
then $X$ has semistable reduction above $v$ over a quadratic 
extension of $F$ if and only if 
either $X$ has an $\I$-invariant point of order $4$,
or all the points of order $2$ on $X$ are $\I$-invariant.
\item[(d)] 
If $p \ne 2$ and $X$ has bad but potentially good reduction
at $v$,
then $X$ has good reduction above $v$ over a quadratic 
extension of $F$ if and only if 
$X$ has no $\I$-invariant point of order $4$ 
and all its points of order $2$ are $\I$-invariant.
\item[(e)] Suppose $p$ is not $2$ or $3$. Then the following
are equivalent:
\begin{enumerate} 
\item[(i)]  $X$ has no $\I$-invariant points of order $2$ or $3$, 
\item[(ii)] there does not exist a finite separable extension $L$
of $F$ of degree less than $6$ 
such that $X$ has semistable reduction at the restriction of 
${\bar v}$ to $L$.
\end{enumerate}
\item[(f)] Suppose $p$ is not $2$ or $3$. Then the following
are equivalent:
\begin{enumerate} 
\item[(i)]  $X$ has no $\I$-invariant points of order $4$ or $3$
and not all the points of order $2$ are $\I$-invariant, 
\item[(ii)] there does not exist a finite separable extension $L$
of $F$ of degree less than $4$ 
such that $X$ has semistable reduction at the restriction 
of ${\bar v}$ to $L$.
\end{enumerate}
\end{enumerate}
\end{cor}

\begin{proof}
Theorem \ref{ellcor} implies that, for $n = 2$, $3$, or $4$, 
if $L$ is an extension of $F$ of degree $R(n)$ which 
is totally ramified above $v$, then 
$X$ has semistable reduction over $L$ above $v$ if and only if 
there exists a subgroup $\s$ of $X_n$ such 
that $\I$ acts as the identity on $\s$ and on $\s^{\perp_n}$.
Parts (a), (b), and (c) are a reformulation of this.

For (d), note that by Theorem 7.4 of \cite{semistab}, if $X$
has an $\I$-invariant point of order $4$ then $X$ has good
reduction at $v$.

In case (e),
if $X$ has an $\I$-invariant point of order $2$ (respectively, $3$),
then $X$ has semistable reduction above $v$ over a totally ramified  
extension of degree $4$ (respectively, $3$), by 
part (a) (respectively, (b)). 
Conversely, suppose $L/F$ is a finite separable extension 
of degree less than $6$, and suppose
$X$ has semistable reduction at the restriction $w$ of ${\bar 
v}$ to $L$.
If $X$ has semistable reduction at $v$, then we are done by 
Corollary \ref{fromdpp} with $n = 6$.
Otherwise, taking completions we have
$[L_{w}:F_{v}] = 2$, $3$, or $4$ 
by Remark \ref{divby23}.
There exists an intermediate
unramified extension $M/F_{v}$ such that $L_{w}/M$ is totally
ramified. 
By parts (a), (b), and (c) applied to $M$ in place of $F$, then 
$X$ has an $\I$-invariant point of order $2$ or $3$.
Case (f) proceeds the same way as case (e). 
\end{proof}

\begin{rem}
Note that if the elliptic curve $X$ has additive reduction at $v$,
but has multiplicative reduction over an extension $L$ of $F$ which
is totally and tamely ramified above $v$, then $X$ has
multiplicative reduction over a quadratic extension of $F$, but
not over any non-trivial 
totally and tamely ramified extension of $F$ of odd degree  
(since $(x+1)^2$ is the only possibility for the characteristic
polynomial of $\rho_\ell(\tau)$, where $\tau$ is as before).
Therefore in case (b) of Corollary \ref{4326cor}, 
either $X$ already has semistable reduction at $v$, 
or else $X$ has good 
(i.e., does not have multiplicative) 
reduction above $v$ over a cubic extension of $F$.
In case (e), $X$ has good reduction
over an extension of degree $6$ or $12$ (see Proposition 1
of \cite{Kraus}).
\end{rem}

\begin{cor}
\label{paddcor}
Suppose 
$X$ is an abelian variety over a field $F$,  
$v$ is a discrete valuation on $F$ of residue characteristic 
$p \ge 0$, and $X$ has purely additive and  
potentially good reduction at $v$. 
\begin{enumerate}
\item[(a)] 
If $p \ne 2$, then $X$ has good reduction above $v$ over a quadratic 
extension of $F$ if and only if 
there exists a subgroup $S$ of $X_4$ such 
that $\I$ acts as the identity on $S$ and on $S^{\perp_4}$.
\item[(b)] If $p \ne 3$, 
then $X$ has good reduction above $v$ over a totally ramified cubic
extension of $F$ if and only if 
there exists a subgroup $S$ of $X_3$ such 
that $\I$ acts as the identity on $S$ and on $S^{\perp_3}$.
\item[(c)] Suppose $p \ne 2$, and
$L/F$ is a degree $4$ extension, 
totally ramified above $v$, 
which has a quadratic subextension over which $X$ has purely
additive reduction.  
Then $X$ has good reduction above $v$ over $L$ if and only if 
there exists a subgroup $S$ of $X_2$ such 
that $\I$ acts as the identity on $S$ and on $S^{\perp_2}$.
\end{enumerate}
\end{cor}

\begin{proof}
The backwards implications follow immediately from 
Corollary \ref{bothcor}.

Let $n = 4$, $3$, and $2$ and $\ell = 2$, $3$, and $2$, 
in cases (a), (b), and (c), respectively. 
Let $\tau$ be a lift to $\I$ of a topological 
generator of the pro-cyclic group $\I/\J$, and let 
$\gamma = \rho_{\ell}(\tau)$.
If $X$ acquires good reduction over a totally ramified degree $R(n)$
extension, then $\gamma^{R(n)} = 1$, by Remark \ref{cyclicrem}. 
Since $X$ has purely additive reduction
at $v$, $1$ is not an eigenvalue of $\gamma$ (see \cite{LenstraOort}).
In case (c), $-1$ is not an eigenvalue of $\gamma$, since $X$
has purely additive reduction over a ramified quadratic extension.
It follows that in cases (a), (b), and (c), respectively,
we have
$$\gamma+1 = 0, \qquad \gamma^2+\gamma+1=0, \qquad 
{\text{ and }}\qquad \gamma^2+1=0$$
in $\End(V_\ell(X))$.
We deduce that $(\gamma-1)^2T_{\ell}(X) \subseteq nT_{\ell}(X)$, i.e.,
$(\tau-1)^2X_n = 0$.
The result now follows from Lemma \ref{sslem}.
\end{proof}

\section{Examples}
\label{exssect}

We will show that the numbers in Corollary \ref{bothcor}
and Theorem \ref{ellcor} are sharp. 

First, we will show that Corollary \ref{bothcor} is sharp in the case of
potentially good reduction. This will show that we cannot take
$Y = X$ in general.
In the next 3 examples, we have $n = 2$, $3$, or $4$, respectively.
Let $\ell$ denote the prime divisor of $n$. 
Suppose that $F$ is a field with a discrete valuation $v$ of 
residue characteristic not equal to $\ell$. Suppose $E$
and $E'$ are elliptic curves over $F$, $E$ has good reduction
at $v$, and $E'$ has additive reduction at $v$ but
acquires good reduction over an extension $L$ of $F$ of degree
$R(n)$.
Let $Y = E \times E'$. 
As shown in the proof of Theorem \ref{bothways}, 
the action of $\I$ on $Y_n$ factors through
$\I/\J$. Let $\tau$ be a lift to $\I$ of a topological 
generator of the pro-cyclic group $\I/\J$, and 
let $g = \rho_{\ell,Y}(\tau)$. Note that $g^{R(n)} = 1$.
Let $G$ denote the cyclic group generated by $g$.
In each example we will construct 
a certain $\Z_{\ell}[G]$-module $T$ such that 
$T \subset T_\ell(Y) \subset \frac{1}{\ell}T$.
Let $C' = \frac{1}{\ell}T/T_\ell(Y)$, 
view $C'$ as a subgroup of $Y_\ell$,
and let $X = Y/C'$. Then $T_\ell(X) \cong T$.
Viewing $T_{\ell}(Y)/T$ as a subgroup
$C$ of $X_{\ell}$, we have $Y = X/C$. 
In our 3 examples, $C$ is stable under $\I$, 
$(\tau - 1)^2X_n \ne 0$, and $(\tau - 1)^2Y_n = 0$. 
By Remark \ref{lemrem}, there is a subgroup $S \subseteq Y_n$
such that $\I$ acts as the identity on $S$ and on $S^{\perp_n}$, 
but there does not exist a subgroup $\s \subseteq X_n$
such that $\I$ acts as the identity on $\s$ and on $\s^{\perp_n}$.
We see that $X$ and $Y$ satisfy (ii) of Corollary \ref{bothcor}.

\begin{ex}
\label{exfor2}
Let $n = 2$. Suppose that $E'$ does not acquire good reduction
over a quadratic subextension of $L/F$. 
As $\Z_{2}[G]$-modules, we have 
$$T_2(Y) \cong 
(\Z_2[x]/(x-1))^2 \oplus \Z_2[x]/(x^2+1),$$
where $g$ acts via multiplication by $x$.
Let 
$$T = \Z_2[x]/(x-1) \oplus \Z_2[x]/(x-1)(x^2+1),$$ and 
view $T$ as a submodule of $T_2(Y)$ via the natural injection.
For example, one could take $F = \Q$, $v = 3$, 
and $E$ and $E'$, respectively, the elliptic curves 11A3 and
36A1 from the tables in \cite{Cremona}.
\end{ex}

\begin{ex}
\label{exfor3}
Let $n = 3$. 
As $\Z_{3}[G]$-modules, we have 
$$T_3(Y) \cong
(\Z_3[x]/(x-1))^2 \oplus \Z_3[x]/(x^2+x+1),$$
where $g$ acts via multiplication by $x$.
Let 
$$T = \Z_3[x]/(x-1) \oplus \Z_3[x]/(x^3-1),$$
and 
view $T$ as a submodule of $T_3(Y)$ via the natural injection.
For example, one could take $F = \Q$, $v = 2$, 
and $E$ and $E'$, respectively, the elliptic curves 11A3 and
20A2 from the tables in \cite{Cremona}.
\end{ex}

\begin{ex}
Let $n = 4$. 
As $\Z_{2}[G]$-modules, we have 
$$T_2(Y) \cong (\Z_2[x]/(x-1))^2 \oplus (\Z_2[x]/(x+1))^2 \cong 
(\Z_2[G])^2,$$
where $g$ acts via multiplication by $x$.
Let 
$$T = \Z_2[x]/(x-1) \oplus \Z_2[x]/(x^2-1) \oplus \Z_2[x]/(x+1),$$ 
and 
view $T$ as a submodule of $T_2(Y)$ via the natural injection.
One could take $F = \Q$, $v = 3$, 
and $E$ and $E'$, respectively, the elliptic curves 11A3 and
99D1 from the tables in \cite{Cremona}.
\end{ex}

Next, we will show that the numbers $8$, $9$ and $4$ (respectively) in
Corollary \ref{bothcor} are sharp.

\begin{ex}
Let $n = 2$, $3$, or $4$. For ease of notation, let $R = R(n)$. 
Let $\ell$ be the prime divisor of $n$. 
Let $F$ be a field with a discrete valuation $v$ of residue 
characteristic not equal to $\ell$, and suppose $E$ is an
elliptic curve over $F$ with multiplicative reduction at $v$.
Suppose that $M$ is a degree $R$ Galois extension of $F$ which is
totally ramified above $v$. Let $\chi$ be the composition
$$\Gal(F^{s}/F) \to \Gal(M/F) \cong \Z/R\Z \hookrightarrow
\Aut_{F}(E^{R}),$$
where the image of the last map is generated by a cyclic permutation
of the factors of $E^{R}$, and $E^{R}$ is the $R$-fold product
of $E$ with itself. 
Let $A$ denote the twist of $E^{R}$ by $\chi$. 
Let $\tau$ denote a lift to $\I$ of a generator of $\I/\J$.
As $\Q_{\ell}[\tau]$-modules, 
$V_{\ell}(A) \cong \Q_{\ell}[\tau]/(\tau^{R}-1)^{2}$.
Let ${\tilde T}$ be the 
inverse image of $\Z_{\ell}[\tau]/(\tau^{R}-1)^{2}$ in $V_{\ell}(A)$.
Then for some integer $k$, we have
$T_{\ell}(A) \subseteq \ell^{k}{\tilde T}$. View 
$\ell^{k}{\tilde T}/T_{\ell}(A)$ as a finite subgroup of $A$
and let $X$ be the quotient of $A$ by this subgroup. 
Then $X$ is defined over an extension $K$ of $F$ unramified above $v$,
and $X$ acquires semistable reduction over $KM$ above $v$. 
We have 
${\tilde T} = T_{\ell}(X)$, and the minimal polynomial of
$\tau$ on $X_{\ell}$ is $(x^{R}-1)^{2} \equiv (x-1)^{2R} \pmod{\ell}$. 
Therefore, 
$$(\tau - 1)^{6}X_{2} \ne 0 \text{ if } n = 2, \,\,
(\tau - 1)^{4}X_{3} \ne 0 \text{ if } n = 3, \text{ and }
(\tau - 1)^{2}X_{2} \ne 0 \text{ if } n = 4.$$   From
 Lemma \ref{lemT} (with $t = 0$,
$F = K$, and $L = KM$) we obtain a lattice $T$ such that 
$$8T \subseteq T_{2}(X) \subseteq T \,\, \text{ if } \,\,  n = 2, \qquad
9T \subseteq T_{3}(X) \subseteq T \,\, \text{ if }\,\,  n = 3,$$ 
$$\text{ and }
\quad 4T \subseteq T_{2}(X) \subseteq T \,\,  \text{ if } \,\,  n = 4.$$
Let $C = T/T_{\ell}(X)$, view $C$ as a subgroup of $X_{\ell}$, 
and let $Y = X/C$. As we saw in the proof of Theorem \ref{bothways},
$(\tau - 1)^{2}Y_{n} = 0$, and
$C$ is killed by $8$, $9$, or $4$ if $n = 2$, $3$, or $4$ respectively.
By Lemma \ref{lemT}efg, the group $C$ is not killed by
$4$, $3$, or $2$, respectively.

Suppose $K'$ is a finite
extension of $K$ unramified above $v$,
$Y'$ is an abelian variety over $K'$, $\varphi : X \to Y'$
is a separable $K'$-isogeny,
and $(\tau - 1)^{2}Y'_{n} = 0$. Suppose that the kernel of
$\varphi$ is killed by some positive integer $s$. Then we can
suppose 
$sT_{\ell}(Y') \subseteq T_{\ell}(X) \subseteq T_{\ell}(Y')$.
 Let $\lambda = (\tau^{2} - 1)/n$.
Since $T_{\ell}(Y')$ is a $\lambda$-stable $\Z_{\ell}$-lattice
in $V_{\ell}(X)$ which contains $T_{\ell}(X)$, we have 
$T \subseteq T_{\ell}(Y')$ by Lemma \ref{lemT}a. 
Therefore, $sT \subseteq T_{\ell}(X)$. 
Then $C$ is killed by $s$, and therefore $s$
cannot be $4$, $3$, or $2$, respectively. This shows
that the numbers $8$, $9$, and $4$ are sharp in Corollary \ref{bothcor}.
Note that $\dim(X) = 4$, $3$, or $2$, respectively. By 
Theorem \ref{ellcor}, these are the smallest dimensions for which
such examples exist.
\end{ex}

\begin{ex}
Let $F$ be a field with a discrete valuation $v$ of residue 
characteristic not equal to $2$, and suppose $E$ is an
elliptic curve over $F$ with multiplicative reduction at $v$.
Suppose that $M$ is a degree $4$ Galois extension of $F$ which is
totally ramified above $v$. Let $\chi$ be the composition
$$\Gal(F^{s}/F) \to \Gal(M/F) \cong \Z/4\Z \hookrightarrow
\Aut_{F}(E^{4}),$$
where the image of the last map is generated by a cyclic permutation
of the factors of $E^{4}$. 
Let 
$$B =
\{(e_{1},e_{2},e_{3},e_{4}) \in E^{4} : 
e_{1} + e_{2} + e_{3} + e_{4} = 0\} \cong E^{3},$$
and let $A$ be the twist of $B$ by $\chi$.
Let $\tau$ denote a lift to $\I$ of a generator of $\I/\J$,
and let $f(x) = (x^{3}+x^{2}+x+1)^{2}$.
As $\Q_{2}[\tau]$-modules, 
$V_{2}(A) \cong \Q_{2}[\tau]/f(\tau)$.
Let ${\tilde T}$ be the 
inverse image of $\Z_{2}[\tau]/f(\tau)$ in $V_{2}(A)$.
As in the previous example, we obtain an abelian variety $X$
such that ${\tilde T} = T_{2}(X)$, and 
such that the minimal polynomial of
$\tau$ on $X_{2}$ is $f(x) \equiv (x-1)^{6}$ (mod $2$). 
Therefore, 
$(\tau - 1)^{4}X_{2} \ne 0$.
As above, we see that $X$ is isogenous 
over an unramified extension to an abelian variety $Y$  
such that $(\tau-1)^{2}Y_{2} = 0$ and such that 
the kernel of the isogeny is killed by $4$.
Using Lemma \ref{lemT}e, we see that there does not
exist such a $Y$ where the kernel is killed by $2$.
This shows that the result in Theorem \ref{ellcor} for
$d=3$ and $n=2$ is sharp. The sharpness of the other numbers 
in Theorem \ref{ellcor} follows from Examples \ref{exfor3} 
and \ref{exfor2}.
\end{ex}

\end{document}